
\input graphicx.tex
\input xy
\xyoption{all}


\magnification\magstephalf

\voffset0truecm
\hoffset=0truecm
\vsize=23truecm
\hsize=15.8truecm
\topskip=1truecm

\binoppenalty=10000
\relpenalty=10000

\font\tenbb=msbm10		\font\sevenbb=msbm7		\font\fivebb=msbm5
\font\tensc=cmcsc10		\font\sevensc=cmcsc7 	\font\fivesc=cmcsc5
\font\tensf=cmss10		\font\sevensf=cmss7		\font\fivesf=cmss5
\font\tenfr=eufm10		\font\sevenfr=eufm7		\font\fivefr=eufm5
\font\grecob=cmmib10


\newfam\bbfam	\newfam\scfam	\newfam\frfam	\newfam\sffam

\textfont\bbfam=\tenbb
\scriptfont\bbfam=\sevenbb
\scriptscriptfont\bbfam=\fivebb

\textfont\scfam=\tensc
\scriptfont\scfam=\sevensc
\scriptscriptfont\scfam=\fivesc

\textfont\frfam=\tenfr
\scriptfont\frfam=\sevenfr
\scriptscriptfont\frfam=\fivefr

\textfont\sffam=\tensf
\scriptfont\sffam=\sevensf
\scriptscriptfont\sffam=\fivesf

\def\bb{\fam\bbfam \tenbb} 
\def\sc{\fam\scfam \tensc} 


\font\sezfont=cmbx10 scaled \magstep1
\font\subsectfont=cmbx10 scaled \magstephalf
\font\titfont=cmbx10 scaled \magstep2
\font\autfont=cmcsc10
\font\intfont=cmss10 

\let\no=\noindent
\let\bi=\bigskip
\let\me=\medskip
\let\sm=\smallskip
\let\ce=\centerline

\let\io=\infty


\newcount\sectno\sectno=0
\newcount\subsectno\subsectno=0
\newcount\thmno\thmno=0
\newcount\tagno\tagno=0
\newcount\notitolo\notitolo=0
\newcount\defno\defno=0

\def\sect#1\par{
	\global\advance\sectno by 1 \global\subsectno=0\global\defno=0\global\thmno=0
	\vbox{\vskip.75truecm\advance\hsize by 1mm
	\hbox{\centerline{\sezfont \the\sectno.~~#1}}
	\vskip.25truecm}\nobreak}

\def\subsect#1\par{
	\global\advance\subsectno by 1
	\vbox{\vskip.75truecm\advance\hsize by 1mm
	\line{\subsectfont \the\sectno.\the\subsectno~~#1\hfill}
	\vskip.25truecm}\nobreak}
	
\def\defin#1{\global\advance\defno by 1
	\global\expandafter\edef\csname+#1\endcsname%
    {\number\sectno.\number\defno}
    \no{\bf Definition~\the\sectno.\the\defno.}}

\def\thm#1#2{
	\global\advance\thmno by 1
	\global\expandafter\edef\csname+#1\endcsname%
	{\number\sectno.\number\thmno}
	\no{\bf #2~\the\sectno.\the\thmno.}}

\def\Tag#1{\global\advance\tagno by 1 {(\the\tagno)}
    \global\expandafter\edef\csname+#1\endcsname%
    		{(\number\tagno)}}
\def\tag#1{\leqno\Tag{#1}}

\def\rf#1{\csname+#1\endcsname\relax}

\def\proof{\no{\sl Proof.}\enskip}
\def\qedn{\thinspace\null\nobreak\hfill\hbox{\vbox{\kern-.2pt\hrule height.2pt
        depth.2pt\kern-.2pt\kern-.2pt \hbox to2.5mm{\kern-.2pt\vrule
        width.4pt \kern-.2pt\raise2.5mm\vbox to.2pt{}\lower0pt\vtop
        to.2pt{}\hfil\kern-.2pt \vrule
        width.4pt\kern-.2pt}\kern-.2pt\kern-.2pt\hrule height.2pt
        depth.2pt \kern-.2pt}}\par\medbreak}
    \def\qed{\hfill\qedn}
    
\newif\ifpage\pagefalse
\newif\ifcen\centrue

\headline={
\ifcen\hfil\else
\ifodd\pageno
\global\hoffset=0.5truecm
\else
\global\hoffset=-0.4truecm
\fi\hfil
\fi}

\footline={
	\ifpage
		\hfill\rm\folio\hfill
	\else
		\global\pagetrue\hfill
\fi}

\lccode`\'=`\'

\def\bib#1{\me\item{[#1]\enskip}}


\def\C{{\bb C}} \def\d{{\rm d}}
 
 \def\N{{\bb N}}

\let\de=\partial

\mathchardef\void="083F
\chardef\BL="15
\def\lambdab{\hbox{\grecob\BL}}

\def\diag{\mathop{\rm Diag}\nolimits}

\def\invlim{\mathop{\vtop{\offinterlineskip
\hbox{\rm lim}\kern1pt\hbox{\kern-1.5pt$\longleftarrow$}\kern-3pt}
}\limits}
\def\casi#1{\vcenter{\normalbaselines\mathsurround=0pt
		\ialign{$##\hfil$&\quad##\hfil\crcr#1\crcr}}}


\ce{\titfont Simultaneous linearization of holomorphic germs} 

\ce{\titfont in presence of resonances}

\me\ce{\autfont Jasmin Raissy}
\sm\ce{\intfont Dipartimento di Matematica, Universit\`a di Pisa}

\ce{\intfont Largo Bruno Pontecorvo 5, 56127 Pisa}

\sm\ce{\intfont E-mail: {\tt raissy@mail.dm.unipi.it}}
\bi

{\narrower

{\sc Abstract.} Let~$f_1, \dots, f_m$ be~$m\ge 2$ germs of biholomorphisms of~$\C^n$, fixing the origin, with~$(\d f_1)_O$ diagonalizable and such that~$f_1$ commutes with~$f_h$ for any~$h=2,\dots, m$. We prove that, under certain arithmetic conditions on the eigenvalues of~$(\d f_1)_O$ and some restrictions on their resonances,~$f_1, \dots, f_m$ are simultaneously holomorphically linearizable if and only if there exists a particular complex manifold invariant under~$f_1, \dots, f_m$. 
}
\bi

\sect Introduction 

One of the main questions in the study of local holomorphic dynamics (see [A] and [B] for general surveys on this topic) is when a given germ of biholomorphism~$f$ of~$\C^n$ at a fixed point~$p$, which we may place at the origin~$O$, is {\it holomorphically linearizable}, i.e., there exists a local holomorphic change of coordinates, tangent to the identity, conjugating~$f$ to its linear part. The answer to this question depends on the set of eigenvalues of~$\d f_O$, usually called the {\it spectrum} of~$\d f_O$. In fact, if we denote by~$\lambda_1, \dots, \lambda_n\in \C^*$ the eigenvalues of~$\d f_O$, then it may happen that there exists a multi-index~$k=(k_1, \dots, k_n)\in \N^n$ with~$|k|:=k_1+\cdots+k_n\ge 2$ and such that
$$\lambda^k - \lambda_j:=\lambda_1^{k_1}\cdots\lambda_n^{k_n} - \lambda_j = 0\tag{eqres}$$
for some~$1\le j\le n$; a relation of this kind is called a {\it resonance} of~$f$, and~$k$ is called a {\it resonant multi-index}. A {\it resonant monomial} is a monomial~$z^k= z_1^{k_1}\cdots z_n^{k_n}$ in the~$j$-th coordinate such that~$\lambda^k = \lambda_j$.

One possible generalization of the previous question is to ask when a given set of~$m\ge 2$ germs of biholomorphisms~$f_1, \dots, f_m$ of~$\C^n$ at the same fixed point, which we may place at the origin, are {\it simultaneously holomorphically linearizable}, i.e., there exists a local holomorphic change of coordinates conjugating~$f_h$ to its linear part for each~$h=1, \dots, m$.

In [R] we found, under certain arithmetic conditions on the eigenvalues and some restrictions on the resonances, a necessary and sufficient condition for holomorphic linearization. In this article we shall use that result to find a necessary and sufficient condition for holomorphic simultaneous linearization.

\sm Before stating our result we need the following definitions:

\sm\defin{De1.1pre} Let~$1\le s\le n$. We say that~$\lambdab = (\lambda_1, \dots, \lambda_s, \mu_1, \dots, \mu_r) \in(\C^*)^n$ {\it has only level~$s$ resonances} if there are only two kinds of resonances:
$$\lambdab^k = \lambda_h \iff k\in\widetilde K_1,\leqno{(a)}$$
where
$$\widetilde K_1 =\left \{k\in \N^n \Bigm | |k|\ge 2,\, \sum_{p=1}^s k_p = 1~~\hbox{and}~~\mu_1^{k_{s+1}}\cdots \mu_r^{k_n}=1 \right\};$$
and
$$\lambdab^k = \mu_j \iff k\in\widetilde K_2,\leqno{(b)}$$
where
$$\widetilde K_2 = \left\{k\in \N^n \Bigm | |k|\ge2, k_1=\cdots=k_s=0~\hbox{and}~\exists j \in\{1, \dots, r\}~\hbox{s.t.}~ \mu_1^{k_{s+1}}\cdots \mu_r^{k_n}=\mu_j \right\}.$$

\sm For~$s=n$ having only level~$s$ resonances means that there are no resonances. When~$s<n$, if~$(\lambda_1, \dots,\lambda_s)$ have no resonances, it is easy to verify that~$\lambdab=(\lambda_1, \dots,\lambda_s, 1, \dots, 1)$ has only level~$s$ resonances.

\sm\defin{De1.0} Let~$n\ge2$ and let~$\lambda_1, \dots, \lambda_n\in\C^*$ be not necessarily distinct. For any~$m\ge 2$ put 
$$\widetilde\omega(m) = \min_{2\le |k|\le m\atop k\not\in {\rm Res}_j(\lambda)} \min_{1\le j\le n} |\lambda^k - \lambda_j|,$$
where~${\rm Res}_j(\lambda)$ is the set of multi-indices~$k\in\N^n$, with $|k|\ge 2$, giving a resonance relation for~$\lambda =(\lambda_1, \dots, \lambda_n)$ relative to~$1\le j\le n$, i.e., such that~$\lambda^k-\lambda_j=0$.
We say that~$\lambda$ {\it satisfies the reduced Brjuno condition} if there exists a strictly increasing sequence of integers~$\{p_\nu\}_{\nu_\ge 0}$ with~$p_0=1$ such that
$$\sum_{\nu\ge 0} p_\nu^{-1} \log\widetilde\omega(p_{\nu+1})^{-1}< \io.$$

\sm Note that the reduced Brjuno condition of order~$n$ (i.e., when there are no resonances) is nothing but the usual Brjuno condition introduced in [Br] (see also [M] pp.~25--37 for the one-dimensional case). 

\sm\defin{De1.1} Let~$f$ be a germ of biholomorphism of~$\C^n$ fixing the origin~$O$ and let~$s\in\N$, with~$1\le s\le n$. The origin~$O$ is called a {\it quasi-Brjuno fixed point of order~$s$} if~$\d f_O$ is diagonalizable and, denoting by~$\lambdab$ the spectrum of~$\d f_O$, we have: {\parindent=30pt
\sm\item{(i)} $\lambdab$ has only level~$s$ resonances;
\sm\item{(ii)} $\lambdab$ satisfies the reduced Brjuno condition.
\sm}
\no We say that~$f$ has the origin as a {\it quasi-Brjuno fixed point} if there exists~$1\le s\le n$ such that it is a quasi-Brjuno fixed point of order~$s$.

\me\defin{DeSimOsc} Let~$f_1, \dots, f_m$ be~$m$ germs of biholomorphisms of~$\C^n$, fixing the origin, with~$m\ge 2$, and let~$M$ be a germ of complex manifold at~$O$ of codimension~$1\le s\le n$, and~$f_h$-invariant for each~$h=1, \dots, m$. We say that~$M$ is a {\it simultaneous osculating manifold for~$f_1, \dots, f_m$} if there exists a holomorphic flat~$(1,0)$-connection~$\nabla$ of the normal bundle~$N_M$ of~$M$ in~$\C^n$ commuting with~$\d f_h|_{N_M}$ for each~$h=1, \dots, m$.

\me In [R] we saw that the osculating condition was necessary and sufficient to extend a holomorphic linearization from an invariant submanifold to a whole neighbourhood of the origin for a germ~$f_1$ of biholomorphism with a quasi-Brjuno fixed point. Our main theorem shows that the simultaneous osculating condition is also necessary and sufficient to extend a common holomorphic linearization, just assuming that~$f_1$ has a quasi-Brjuno fixed point and commutes with~$f_2,\dots, f_m$:

\sm\thm{TeoremaIntro}{Theorem} {\sl Let~$f_1, \dots, f_m$ be~$m\ge 2$ germs of biholomorphisms of~$\C^n$, fixing the origin. Assume that~$f_1$ has the origin as a quasi-Brjuno fixed point of order~$s$, with~$1\le s\le n$, and that it commutes with~$f_h$ for any~$h=2,\dots, m$. Then~$f_1, \dots, f_m$ are simultaneously holomorphically linearizable if and only if there exists a germ of complex manifold~$M$ at~$O$ of codimension~$s$, invariant under~$f_h$ for each~$h=1, \dots, m$, which is a simultaneous osculating manifold for~$f_1, \dots, f_m$ and such that~$f_1|_M, \dots, f_m|_M$ are simultaneously holomorphically linearizable.} 

\sm A similar topic is studied in [S]. However, his results are not comparable with ours, because his notion of ``linearization modulo an ideal'' is not suitable for producing a full linearization result, except when there are no resonances at all, whereas in our result we explicitly admit some resonances.

We shall need the following notation: if~$g \colon\C^n\to \C$ is a holomorphic function with~$g(O)=0$, and~$z=(x, y)\in \C^n$ with~$x\in\C^s$ and~$y\in\C^{n-s}$, we shall denote by~${\rm ord}_x(g)$ the maximum positive integer~$m$ such that~$g$ belongs to the ideal~$\langle x_1,\cdots, x_s\rangle^m$. Furthermore, we shall say that the local coordinates~$z=(x,y)$ are {\it adapted} to the complex submanifold~$M$ if in those coordinates~$M$ is given by~$\{x=0\}$.

\sect Linearization

We first introduced {\it osculating manifolds} in [R]. A germ~$f$ of biholomorphism of~$\C^n$ fixing the origin~$O$ admits an {\it osculating manifold~$M$ of codimension~$1\le s\le n$} if there is a germ of~$f$-invariant complex manifold~$M$ at~$O$ of codimension~$s$ such that the normal bundle~$N_M$ of~$M$ admits a holomorphic flat~$(1,0)$-connection that commutes with~$\d f|_{N_M}$. Definition \rf{DeSimOsc} is the natural extension of this object to the case we are dealing with. 

\sm We shall need the following characterization of simultaneous osculating manifolds.

\sm\thm{PrSimOsc}{Proposition} {\sl Let~$f_1, \dots, f_m$ be~$m$ germs of biholomorphisms of~$\C^n$, fixing the origin, with~$m\ge 2$, and let~$M$ be a germ of complex manifold at~$O$ of codimension~$1\le s\le n$, and~$f_h$-invariant for each~$h=1, \dots, m$. Then~$M$ is a simultaneous osculating manifold for~$f_1, \dots, f_m$ if and only if there exist local holomorphic coordinates~$z=(x,y)$ about~$O$ adapted to~$M$ in which~$f_h$ has the form
$$\casi{\displaystyle{x_i'= \sum_{p=1}^s a_{i,p}^{(h)}x_p + \widehat f^{(h)}_i(x,y)} &\hbox{for}~$i=1,\dots, s$,\cr\noalign{\sm}
 			y_j'=f^{(h)}_j(x,y) &\hbox{for}~$j=1,\dots, r=n-s$,}\tag{osc}$$
with
$${\rm ord}_x(\widehat f_i^{(h)})\ge 2,$$
for any~$i=1, \dots, s$ and~$h=1,\dots, m$.}

\sm\proof If there exist local holomorphic coordinates~$z=(x,y)$ about~$O$ adapted to~$M$ in which~$f_h$ has the form \rf{osc} with~${\rm ord}_x(\widehat f_i^{(h)})\ge 2$ for any~$i=1, \dots, s$ and~$h=1,\dots, m$, then it is obvious to verify that the trivial holomorphic flat~$(1,0)$-connection commutes with~$\d f_h|_{N_M}$ for each~$h=1, \dots, m$.

Conversely, let~$\nabla$ be a holomorphic flat~$(1,0)$-connection of the normal bundle~$N_M$ commuting with~$\d f_h|_{N_M}$ for each~$h=1, \dots, m$.
It suffices to choose local holomorphic coordinates~$z=(x,y)$ adapted to~$M$ in which all the connection coefficients~$\Gamma^i_{jk}$ with respect to the local holomorphic frame~$\{\pi({\de\over\de x_1}), \dots, \pi({\de\over\de x_s})\}$ of~$N_M$ are zero (see [R] Proposition~3.1 and Lemma~3.2), and then the assertion follows immediately from the proof of Theorem~$1.3$ of [R].
\qed

\sm\thm{CoSimOsc}{Corollary} {\sl Let~$f_1, \dots, f_m$ be~$m$ germs of biholomorphisms of~$\C^n$, fixing the origin, with~$m\ge 2$, and let~$M$ be a germ of complex manifold at~$O$ of codimension~$1\le s\le n$, and~$f_h$-invariant for each~$h=1, \dots, m$. Then~$M$ is a simultaneous osculating manifold for~$f_1, \dots, f_m$ such that~$f_1|_M, \dots, f_m|_M$ are simultaneously holomorphically linearizable if and only if there exist local holomorphic coordinates~$z=(x,y)$ about~$O$ adapted to~$M$ in which~$f_h$ has the form
$$\casi{\displaystyle{x_i'= \sum_{p=1}^s a_{i,p}^{(h)}x_h + \widehat f^{(h),1}_i(x,y)} &\hbox{for}~$i=1,\dots, s$,\cr\noalign{\sm}
  		\displaystyle{y_j'=f^{(h){\rm lin}}_j(x,y)+ \widehat f^{(h),2}_j(x,y)} &\hbox{for}~$j=1,\dots, r=n-s$,}\tag{osc3}$$
where~$f^{(h){\rm lin}}_j(x,y)$ is linear and
$$\eqalign{&{\rm ord}_x (\widehat f_i^{(h),1})\ge 2,\cr
  			&{\rm ord}_x (\widehat f_j^{(h),2})\ge 1, }\tag{osc2}$$
for any~$i=1, \dots, s$,~$j=1,\dots, r$ and~$h=1,\dots, m$.}

\sm\proof One direction is clear.

Conversely, thanks to Proposition \rf{PrSimOsc}, the fact that~$M$ is a simultaneous osculating manifold for~$f_1, \dots, f_m$ is equivalent to the existence of local holomorphic coordinates~$z=(x,y)$ about~$O$ adapted to~$M$, in which~$f_h$ has the form \rf{osc3}\ with~${\rm ord}_x (\widehat f_i^{(h),1})\ge 2$ for any~$i=1, \dots, s$ and~$h=1,\dots, m$. 
%
Furthermore,~$f_1|_M, \dots, f_m|_M$ are simultaneously holomorphically linearizable; therefore there exists a local holomorphic change of coordinate, tangent to the identity, and of the form
$$\eqalign{&\tilde x= x,\cr
  		   &\tilde y = \Phi(y),}$$
conjugating~$f_h$ to~$\tilde f_h$ of the form \rf{osc3} satisfying \rf{osc2}, for each~$h=1,\dots, m$, as we wanted. \qed

\sm\thm{ReFormal}{Remark} It is possible to give the formal analogous of Definition~\rf{DeSimOsc}, and then to prove a formal analogous of Proposition \rf{PrSimOsc} and Corollary \rf{CoSimOsc}, exactly as in [R]. 

\me In the proof of Theorem \rf{TeoremaIntro} we shall use the following result we proved in [R]

\sm\thm{Teorema1}{Theorem} (Raissy, 2007) {\sl Let~$f$ be a germ of biholomorphism of~$\C^n$ having the origin~$O$ as a quasi-Brjuno fixed point of order~$s$. Then~$f$ is holomorphically linearizable if and only if it admits an osculating manifold~$M$ of codimension~$s$ such that~$f|_M$ is holomorphically linearizable.}

\sm We can now prove our result.

\sm\thm{Teorema}{Theorem} {\sl Let~$f_1, \dots, f_m$ be~$m\ge 2$ germs of biholomorphisms of~$\C^n$, fixing the origin. Assume that~$f_1$ has the origin as a quasi-Brjuno fixed point of order~$s$, with~$1\le s\le n$, and that it commutes with~$f_h$ for any~$h=2,\dots, m$. Then~$f_1, \dots, f_m$ are simultaneously holomorphically linearizable if and only if there exists a germ of complex manifold~$M$ at~$O$ of codimension~$s$, invariant under~$f_h$ for each~$h=1, \dots, m$, which is a simultaneous osculating manifold for~$f_1, \dots, f_m$ and such that~$f_1|_M, \dots, f_m|_M$ are simultaneously holomorphically linearizable.} 

\sm\proof Let~$M$ be a germ of complex manifold at~$O$ of codimension~$s$, invariant under~$f_h$ for each~$h=1, \dots, m$ which is a simultaneous osculating manifold for~$f_1, \dots, f_m$ and such that~$f_1|_M, \dots, f_m|_M$ are simultaneously holomorphically linearizable. Thanks to the hypotheses we can choose local holomorphic coordinates
$$(x,y)=(x_1, \dots, x_s, y_1, \dots, y_r)$$
such that~$f_1$ is of the form
$$\eqalign{&x_i'=\lambda_{1,i} x_i + f^{(1),1}_i(x,y) \quad\hbox{for}~i=1,\dots, s,\cr
  			&y_j'=\mu_{1,j} y_j + f^{(1),2}_j(x,y) \quad\hbox{for}~j=1,\dots, r=n-s,}$$
and, for~$h=2,\dots, m$, each~$f_h$ is of the form
$$\casi{\displaystyle{x_i'=\sum_{p=1}^s a_{i,p}^{(h)}x_p  + f^{(h),1}_i(x,y) }&\hbox{for}~$i=1,\dots, s$,\cr\noalign{\sm}
  			\displaystyle{y_j'=f^{(h){\rm lin}}_j(x,y) + f^{(h),2}_j(x,y)} &\hbox{for}~$j=1,\dots, r=n-s$,}$$
where~$f^{(h){\rm lin}}_j(x,y)$ is linear, and for each~$k=1,\dots,m$
$$\eqalign{&{\rm ord}_x(f^{(k),1}_i)\ge 2,\cr
  			&{\rm ord}_x(f^{(k),2}_j)\ge 1,}$$
that is
$$\eqalign{&f^{(k),1}_i(x,y)=\sum_{|K|\ge 2\atop |K'|\ge 2} f^{(k),1}_{K,i}x^{K'}y^{K''} \quad\hbox{for}~i=1,\dots, s,\cr
  			&f^{(k),2}_j(x,y)=\sum_{|K|\ge 2\atop |K'|\ge 1} f^{(k),2}_{K,j}x^{K'}y^{K''} \quad\hbox{for}~j=1,\dots, r,}$$
where~$K=(K',K'')\in\N^s\times\N^{r}=\N^n$ and~$|K|=\sum_{p=1}^n K_p$.

Thanks to Theorem \rf{Teorema1} and its proof, we know that~$f_1$ is holomorphically linearizable via a linearization~$\psi$ of the form
$$\eqalign{&x_i= u_i + \psi^1_i(u,v) \quad\hbox{for}~i=1,\dots, s,\cr
  			&y_j= v_j + \psi^2_j(u,v) \quad\hbox{for}~j=1,\dots, r, }$$
where~$(u,v)=(u_1, \dots, u_s,v_1, \dots,v_r)$ and
$$\eqalign{&{\rm ord}_u(\psi_i^1)\ge 2,\cr
  			&{\rm ord}_u(\psi_j^2)\ge 1,}$$
that is
$$\eqalign{&\psi^{1}_i(u,v)=\sum_{|K|\ge 2\atop |K'|\ge 2} \psi^{1}_{K,i}u^{K'}v^{K''} \quad\hbox{for}~i=1,\dots, s,\cr
  			&\psi^{2}_j(u,v)=\sum_{|K|\ge 2\atop |K'|\ge 1} \psi^{2}_{K,j}u^{K'}v^{K''} \quad\hbox{for}~j=1,\dots, r.}$$
Since~$\psi^{-1}\circ f_1\circ\psi=\diag(\lambda_{1,1},\dots, \lambda_{1,s},\mu_{1,1},\dots, \mu_{1,r})$ commutes with~$\tilde f_h= \psi^{-1}\circ f_h\circ\psi$ for each~$h=2, \dots, m$, and~$(\lambda_{1,1},\dots, \lambda_{1,s},\mu_{1,1},\dots, \mu_{1,r})$ has only level~$s$ resonances, it is immediate to verify that~$\tilde f_h$ has the form
$$\casi{\displaystyle{u_i'=\sum_{p=1}^s a_{i,p}^{(h)}u_p  + \sum_{1\le l\le n\atop \lambda_{1,l}=\lambda_{1,i}}u_l\tilde f^{(h),1}_{l,i}(v)} &\hbox{for}~$i=1,\dots, s$,\cr\noalign{\sm}
  			\displaystyle{v_j'= f^{(h){\rm lin}}_j(u,v) + \tilde f^{(h),2}_j(v)} &\hbox{for}~$j=1,\dots, r$.}$$
Moreover, since~$f_h \circ\psi = \psi\circ\tilde f_h$, we have
$$\eqalign{&\!\!\!{\sum_{p=1}^s a^{(h)}_{i,p}\sum_{|K|\ge 2\atop |K'|\ge 2} \psi^{1}_{K,p}u^{K'}v^{K''}+ \sum_{|K|\ge 2\atop |K'|\ge 2} f^{(h),1}_{K,i}(u+\psi^1(u,v))^{K'}(v+\psi^2(u,v))^{K''}} \cr 
&\qquad= \sum_{1\le l\le n\atop \lambda_{1,l}=\lambda_{1,i}}\!\!\!\!u_l\tilde f^{(h),1}_{l,i}(v) \cr
&\quad\qquad+ \!\!\!\sum_{|K|\ge 2\atop |K'|\ge 2} \!\!\psi^{1}_{K,i}\left(\!\sum_{p=1}^s a_{1,p}^{(h)}u_p + \!\!\!\!\sum_{1\le l\le n\atop \lambda_{1,l}=\lambda_{1,1}}\!\!\!\!u_l\tilde f^{(h),1}_{l,1}(v)\!\right)^{K_1}\!\!\!\!\!\!\!\cdots \left(\!\sum_{p=1}^s a_{s,p}^{(h)}u_p  + \!\!\!\!\sum_{1\le l\le n\atop \lambda_{1,l}=\lambda_{1,s}}\!\!\!\!u_l\tilde f^{(h),1}_{l,s}(v)\!\right)^{K_s} \cr
&\qquad\qquad\qquad\qquad\qquad\qquad\qquad\qquad\qquad\quad\quad\quad\quad\quad\quad\quad\times(f^{(h){\rm lin}}(u,v)+ \tilde f^{(h),2}(v))^{K''}}\tag{eqlin1}$$
for~$i=1,\dots, s$, and
$$\eqalign{&\!\!\!\sum_{q=1}^r b^{(h)}_{j,q}\sum_{|K|\ge 2\atop |K'|\ge 1} \psi^{2}_{K,q}u^{K'}v^{K''}+ \sum_{p=1}^s c^{(h)}_{j,p}\sum_{|K|\ge 2\atop |K'|\ge 2} \psi^{1}_{K,p}u^{K'}v^{K''}\cr
&\qquad\qquad\qquad\qquad\qquad\qquad\qquad\qquad+ \sum_{|K|\ge 2\atop |K'|\ge 1} f^{(h),2}_{K,j}(u+\psi^1(u,v))^{K'}(v+\psi^2(u,v))^{K''} \cr 
&\qquad= \tilde f^{(h),2}_j(v)\cr
&\quad\qquad+ \!\!\!\sum_{|K|\ge 2\atop |K'|\ge 1} \!\!\psi^{2}_{K,i}\left(\!\sum_{p=1}^s a_{1,p}^{(h)}u_p + \!\!\!\!\sum_{1\le l\le n\atop \lambda_{1,l}=\lambda_{1,1}}\!\!\!\!u_l\tilde f^{(h),1}_{l,1}(v)\!\right)^{K_1}\!\!\!\!\!\!\!\cdots \left(\!\sum_{p=1}^s a_{s,p}^{(h)}u_p  + \!\!\!\!\sum_{1\le l\le n\atop \lambda_{1,l}=\lambda_{1,s}}\!\!\!\!u_l\tilde f^{(h),1}_{l,s}(v)\!\right)^{K_s} \cr
&\qquad\qquad\qquad\qquad\qquad\qquad\qquad\qquad\qquad\quad\quad\quad\quad\quad\quad\quad\times(f^{(h){\rm lin}}_j(u,v)+ \tilde f^{(h),2}(v))^{K''}}\tag{eqlin2}$$
for~$j=1,\dots, r$. 

Now, it is not difficult to verify that there are no terms of the form~$u^{K'}v^{K''}$ with~$|K'|=1$ in the left-hand side of \rf{eqlin1}, whereas in the right-hand side terms of this form are given only by the sum of the~$u_l\tilde f^{(h),1}_{l,i}(v)$; therefore it must be
$$\tilde f^{(h),1}_{l,i}(v)\equiv 0,$$
for all pairs $l$, $i$. Similarly, there are no terms of the form~$u^{K'}v^{K''}$ with~$K'=O$ in the left-hand side of \rf{eqlin2}, whereas, again, in the right-hand terms of this form are given by~$\tilde f^{(h),2}_j(v)$ only; so
$$\tilde f^{(h),2}_j(v)\equiv 0\quad\hbox{for}~j=1,\dots, r.$$
This proves that~$\tilde f_h$ is linear for every~$h=2,\dots, m$, that is~$\psi$ is a simultaneous holomorphic linearization for~$f_1, \dots, f_m$. 

The other direction is clear. In fact, if~$f_1$ commutes with~$f_2,\dots, f_m$  and~$f_1,\dots,f_m$ are linear, then the eigenspace of~$f_1$ relative to the eigenvalues~$\mu_{1,1}, \dots, \mu_{1,r}$ is a simultaneous osculating manifold for~$f_1, \dots, f_m$ (and~$f_1|_M, \dots, f_m|_M$ are linear), where~$(\lambda_{1,1}, \dots, \lambda_{1,s}, \mu_{1,1}, \dots, \mu_{1,r})$ is the spectrum of~$f_1$. \qed

\sm\thm{CoSimLin}{Corollary} {\sl Let~$f_1, \dots, f_m$ be~$m\ge 2$ germs of commuting biholomorphisms of~$\C^n$, fixing the origin. Assume that~$f_1$ has the origin as a quasi-Brjuno fixed point of order~$s$, with~$1\le s\le n$. Then~$f_1, \dots, f_m$ are simultaneously holomorphically linearizable if and only if there exists a germ of complex manifold~$M$ at~$O$ of codimension~$s$, invariant under~$f_h$ for each~$h=1, \dots, m$ which is a simultaneous osculating manifold for~$f_1, \dots, f_m$ and such that~$f_1|_M, \dots, f_m|_M$ are simultaneously holomorphically linearizable.} 


\sm As a final corollary, taking~$s=n$ in Theorem \rf{Teorema}, one gets

\sm\thm{CoBrjuno}{Corollary} {\sl Let~$f_1, \dots, f_m$ be~$m\ge 2$ germs of biholomorphisms of~$\C^n$, fixing the origin. Assume that~$f_1$ has the origin as a Brjuno fixed point, and that it commutes with~$f_h$ for any~$h=2,\dots, m$. Then~$f_1, \dots, f_m$ are simultaneously holomorphically linearizable.}

\vbox{\vskip.75truecm\advance\hsize by 1mm
	\hbox{\centerline{\sezfont References}}
	\vskip.25truecm}\nobreak
\parindent=40pt

\bib{A} {\sc Abate, M.:} {\sl Discrete holomorphic local dynamical systems,} to appear in ``Holomorphic Dynamical Systems'', Eds. G. Gentili, J. Guenot, G. Patrizio, Lectures notes in Math., Springer-Verlag, Berlin, 2009, arXiv:0903.3289v1.

\bib{B} {\sc Bracci, F.:} {\sl Local dynamics of holomorphic diffeomorphisms,} Boll. UMI (8), 7--B (2004), pp. 609--636.

\bib{Br} {\sc Brjuno, A. D.:} {\sl Analytic form of differential equations,} Trans. Moscow Math. Soc., {\bf 25} (1971), pp. 131--288; {\bf 26} (1972), pp. 199--239.

\bib{M} {\sc Marmi, S.:} ``An introduction to small divisors problems'', I.E.P.I., Pisa, 2003.

\bib{R} {\sc Raissy, J.:} {\sl Linearization of holomorphic germs with quasi-Brjuno fixed points,}  Math. Z. (2009), {\intfont http://www.springerlink.com/content/3853667627008057/fulltext.pdf}, On\-li\-ne First.

\bib{S} {\sc Stolovitch, L.:} {\sl Family of intersecting totally real manifolds of $(\C^n,0)$ and CR-singularities,} preprint 2005, (arXiv: math/0506052v2).

\bye